\documentclass[11pt]{amsart}
\usepackage{ amsmath}
\usepackage{amssymb}
\usepackage{amsxtra}
\usepackage{enumerate}
\usepackage[mathscr]{eucal}
\usepackage[margin=3.3 cm]{geometry}

\def\R{\mathbb R}

\def\H{\mathbb H}
\def \h{\bf H}

\def\SL{{\rm SL}}
\def\at{{\rm argtr}}
\def\abt{{\rm abstr }}
\def\C{\mathbb C}

\def\S{{\rm SL}(2, \H)}

\def \s{\mathbb S}
\newtheorem{theorem}{Theorem}[section]
\newtheorem{lemma}[theorem]{Lemma}

\theoremstyle{definition}
\newtheorem{definition}{Definition}

\theoremstyle{remark}

\numberwithin{equation}{section}
\theoremstyle{plain}

\newtheorem{cor}[theorem]{Corollary}

\newtheorem*{ack}{Acknowledgement}

\newcommand{\secref}[1]{Section~\ref{#1}}
\newcommand{\thmref}[1]{Theorem~\ref{#1}}
\newcommand{\lemref}[1]{Lemma~\ref{#1}}

\newcommand{\corref}[1]{Corollary~\ref{#1}}
\newcommand{\eqnref}[1]{~{\textrm(\ref{#1})}}

\begin{document}
\title[ Test Map and Discreteness in $\S$]{ Test map and Discreteness in $\S$}
\author[K. Gongopadhyay]{Krishnendu Gongopadhyay}
 \address{ Indian Institute of Science Education and Research (IISER) Mohali,
Knowledge City, Sector 81, SAS Nagar, Punjab 140306, India}
\email{krishnendu@iisermohali.ac.in, krishnendug@gmail.com}

\author[A. Mukherjee]{Abhishek Mukherjee}
\address{ Kalna College, Kalna,  Dist. Burdwan, West Bengal 713409, and}
\address{Department of Mathematics, Jadavpur University, Jadavpur, Kolkata 700032 }
\email{abhimukherjee.math10@gmail.com }

\author[S. K. Sardar]{Sujit Kumar Sardar}
\address{Department of Mathematics, Jadavpur University, Jadavpur, Kolkata 700032 }
\email{sksardar@math.jdvu.ac.in }


\date{\today}

%\tableofcontents
\begin{abstract}
Let $\H$ be the division ring of real quaternions. 
Let $\SL(2,\H)$ be the group of  $2 \times 2$ quaternionic matrices $A=\begin{pmatrix} a & b \\ c & d \end{pmatrix}$ with quaternionic determinant $\det A=|ad-aca^{-1} b|=1$. This group acts by the  orientation-preserving isometries of the five dimensional real hyperbolic space. We obtain discreteness criteria for Zariski-dense subgroups of $\SL(2, \H)$. 
\end{abstract}
\maketitle

\section{Introduction}

  Let $\h^{n+1}$ be the $(n+1)$-dimensional (real) hyperbolic space and let $M(n)$ denote the (orientation-preserving) M\"obius group that acts on $\h^{n+1}$ by isometries. Given a subgroup $G$ of $M(n)$, it is an interesting problem to ask when $G$ is discrete. In particular, one asks when a two-generator subgroup of $M(n)$ is discrete. It has been seen in the literature, especially for $n=2$,  that the discreteness of the two-generator subgroups of $G$ determine the discreteness of $G$. The linear group $\SL(2, \C)$ acts on $\partial \h^3\approx \s^2$ by linear fractional transformations,  and this action identifies the group $M(2)$ with  ${\mathrm{PSL}}(2, \C)$,   e.g. see \cite{beardon}. The J{\o}rgensen inequality in $\SL(2, \C)$ gave a sufficient algorithm for discreteness of a two-generator subgroup. There have been many attempts in the literature to formulate generalizations of J{\o}rgensen inequality in higher dimensions and to obtain discreteness criteria using two-generator subgroups,  e.g. see  \cite{fn}, \cite{lw},  \cite{TW02}, \cite{wlc}, \cite{yang0} for some recent investigations in this direction.

 A subgroup $G$ of $M(n)$ is called \emph{Zariski-dense} if it does not have a global fixed point and neither it preserves a proper totally geodesic subspace of $\h^{n+1}$. In \cite{ah}  Abikoff and Haas proved that a Zariski-dense subgroup $G$ of $M(n)$ is discrete if and only if every two-generator subgroup $\langle f, g \rangle$ of $G$ is discrete. When $n$ even,  Abikoff and Haas proved a stronger result that says that a Zariski-dense subgroup $G$ of $M(2m)$ is discrete if and only if every cyclic subgroup of $G$ is discrete. This implies that the discreteness of a subgroup in $M(2m)$ is controlled by the cyclic subgroups. In \cite{chen}, Chen obtained a discreteness criterion that uses a fixed (test) map to check  discreteness of a M\"obius subgroup. Chen proved that a Zariski-dense subgroup $G$ of $M(n)$ is discrete if for any $g$ in $G$, and a fixed non-trivial element $f$ from $M(n)$, the group $\langle f, g \rangle$ is discrete, where $f$ is  not an irrational rotation (i.e. of infinite order) or if having finite order, it acts as a non-identity M\"obius transformation on the minimal sphere containing the limit set of $G$. Chen's discreteness criterion  involves two-generator subgroups of $M(n)$ with only one generator from $G$ itself. 

Motivated by Chen's work it is natural to ask how far the test map $f$ may be chosen outside $G$.  This was the line of investigation of  Yang who asked this problem for $\SL(2, \C)$ in \cite{yang1}. Yang gave a partial answer to this question  and formulated a conjecture for the remaining cases.  In \cite{cao}, Cao completed Yang's program by solving Yang's conjecture.  Yang and Zhao \cite{yang2} gave another proof to the conjecture.  Recently, Yang and Zhao  \cite{yang3} have obtained a discreteness criterion in $\SL(2, \C)$ that says that a non-elementary subgroup $G$ of $\SL(2, \C)$ is discrete if every two generator subgroup $\langle g, fgf^{-1} \rangle$ is discrete, where $g$ is a non-trivial element of $G$ and $f$ is an arbitrary but fixed  element in $\SL(2, \C)$. The work of Cao and Yang et. al.  shows that the discreteness of a subgroup $G$ of $\SL(2, \C)$ is completely determined by two-generator subgroups $\langle f, g \rangle$, where $f$ is a test map and $g$ is an element of $G$. However, given a test map $f$,  it is not clear from these works that whether the elements $g$ from $G$ can be restricted to a smaller class. 

\medskip The aim of this paper is to investigate the above problems in higher dimensions. We focus on the group $M(4)$ that provides the closest analogue of  ${\rm PSL}(2, \C)$ action on the Riemann sphere by M\"obius transformations.  Let $\H$ be the division ring of real quaternions. 
Let $\SL(2,\H)$ be the group of  $2 \times 2$ quaternionic matrices $A=\begin{pmatrix} a & b \\ c & d \end{pmatrix}$ with quaternionic determinant $\det A=|ad-aca^{-1} b|=1$. The group ${\rm PSL}(2, \H)=\SL(2, \H)/\{ \pm I\}$ can be identified with the group of orientation-preserving isometries of the five dimensional hyperbolic space using the quaternionic linear fractional transformations, see \cite{bg}, \cite{p}, \cite{wi}. We investigate the discreteness of two-generator subgroups using this action. 

\medskip  To state our main results we recall from \cite{kg}, \cite{p} that  a parabolic element in $\SL(2, \H)$  is conjugate to
\begin{equation}\label{1} \begin{pmatrix} \lambda & 1 \\ 0 & \lambda \end{pmatrix}, ~ |\lambda|=1, ~\lambda \in \C\end{equation} 
and upto conjugacy, an elliptic or hyperbolic  element $A$  is  given by  
\begin{equation}\label{2} A=\begin{pmatrix} \lambda & 0 \\ 0 & \mu \end{pmatrix},\end{equation} 
where $\lambda, \mu \in \C$, and $A$ is hyperbolic if and only if  $|\lambda| \neq 1 \neq |\mu|$.  If $|\lambda|=|\mu| =1$ and $\lambda$ is not similar to $\mu$ in $\H^{\ast}$, then $A$ is called \emph{2-rotatory elliptic}. 
\begin{definition}
Let $A$ be an elliptic or hyperbolic element in $\S$ which is represented by \eqnref{1} or \eqnref{2} up to conjugacy.  We define the \emph{argument trace} of $A$ by $${\at(A)}=\arg(\lambda) + \arg(\mu),$$
and the \emph{absolute trace} of $A$ by:
$${\abt(A)}=|\lambda| + |\mu|.$$
\end{definition}
Note that an element of $\S$ is hyperbolic if and only if ${\abt(A)} > 2$. Now we state our main result. 

\begin{theorem}\label{thm1}
Let $G$ be a  Zariski-dense subgroup of $\S$. 
\begin{enumerate}
\item Let $f$ be a $2$-rotatory elliptic element of $\S$ such that  {$0<\at(f) <\frac{\pi}{3}$}.   If  the two generator subgroup $\langle f,g\rangle$ is discrete for  every hyperbolic element $g$ in $G$, then $G$ is discrete.

\medskip \item Let $f$ be a hyperbolic element of $\S$ such that $$ \frac{1}{2}(\abt^2(f) -3) < \cos(\at(f)).$$
If the two generator subgroup $\langle f,g \rangle $ is  discrete for  every hyperbolic element $g$ in $G$, then $G$ is discrete.

\medskip \item    Let $f$ be a parabolic  element of $\S$ such that,  up to conjugacy,
$$f = \begin{pmatrix} 1 & \mu\\ 0& 1\end{pmatrix}, \;  |\mu| \leq 1 .$$
 If the two generator subgroup $\langle f,g \rangle $ is  discrete  for every hyperbolic element $g$ in $G$, then $G$ is discrete.
\end{enumerate}
\end{theorem}
After proving the above result, using similar methods we have obtained the following.
\begin{theorem}\label{thm2}
Let $G$ be a  Zariski-dense subgroup of $\S$.
\begin{enumerate}
\item Let $f$ be a $2$-rotatory elliptic element of $\SL(2, \H)$ such that  {$0<\at(f) <\frac{\pi}{3}$}.   If the two generator subgroup $\langle f,gfg^{-1}\rangle$ is discrete and non-elementary for  every hyperbolic element $g$ in $G$, then $G$ is discrete.
\medskip\item Let $f$ be a  hyperbolic  element of $\SL(2, \H)$ such that
 $$ \frac{1}{2}(\abt^2(f) -3) < \cos(\at(f)).$$  
 If  the two generator subgroup $\langle f,gfg^{-1}\rangle$ is discrete for  every hyperbolic element $g$ in $G$, then $G$ is discrete.
\medskip \item    Let  $f$ be a parabolic  element of $\SL(2, \H)$ such that,  up to conjugacy,
$$f = \begin{pmatrix} 1 & \mu\\ 0& 1\end{pmatrix}, \;  |\mu| \leq 1 .$$ 
 If the two generator subgroup $\langle f,gfg^{-1}\rangle$ is discrete for  every hyperbolic element $g$ in $G$, then $G$ is discrete.
\end{enumerate}
\end{theorem}
The above two  theorems indicate that  the discreteness of a Zariski-dense subgroup $G$ of $\SL(2, \H)$, equivalently, $M(n)$, $n \leq 5$, is determined  by the  two-generator subgroups  involving a test map and the hyperbolic elements of $G$.  It is interesting to note that our choice of $f$ in $\S$ lies in a very nice region where one can choose uncountably many irrational rotations which are of infinite orders. Given the dynamical type of the test map, it belongs to an one parameter family where each element in the family may be chosen as a test map.

We note here that the restrictions on $\at(f)$ and $\abt(f)$ in both the theorems are necessary. These quantities come from the J{\o}rgensen type inequalities in \cite{gm} and cannot be relaxed. In part (1) of both the theorems, the quantity $\at(f)$ can not be zero, as in that case $f$ will reduces to an $1$-rotatory elliptic. If $\at(f)=\frac{\pi}{3}$ , then the arguments we give here become inconclusive.  Similarly in part (2), equality of the given inequality would imply that $f$ is an elliptic of order at least seven, by \cite[Corollary 8]{gm}. This would contradict the hypothesis that $f$ is hyperbolic.

Plan of the paper is as follows. In \secref{prel} we recall some preliminary results that  include  J{\o}rgensen type inequalities for two generator subgroups of $\SL(2, \H)$ as obtained in \cite{gm}, also see \cite{kel}, \cite{wat}. We apply these results to prove \thmref{thm1} and \thmref{thm2} in \secref{ms}.

\section{Preliminaries}\label{prel}
\subsection{The Quaternions}Let $\H$ denote the division ring of quaternions. Recall that every element of $\H$ is of the form  $a_{0}+a_{1}i+a_{2}j+a_{3}k$,where $a_{0},a_{1},a_{2},a_{3}\in \R$, and  $ i,j,k$ satisfy relations:  $i^{2}=j^{2}=k^{2}=-ijk=-1$. Any $a\in {\H}$ can be uniquely written as  $a=a_{0}+a_{1}i+a_{2}j+a_{3}k$.
We define $\Re(a)=a_{0}$=the real part of $a$ and $\Im(a)=a_{1}i+a_{2}j+a_{3}k=$ the imaginary part of $a$. Also, define the conjugate of $a$ as $\overline {a}= \Re(a)-\Im(a)$.
The norm of $a$ is $|a|=\sqrt{a_0^{2}+a_1^{2}+a_2^{2}+a_3^{2}}$.
Two quaternions ${a,b}$ are said to be \emph{similar} if there exists a non-zero quaternion $ {c}$ such that $ {b=c^{-1}ac}$ and we write it as $ {a\backsim b}$. It is easy to verify that $ {a \backsim b}$ if and only if $ {\Re(a)=\Re(b)}$ and $|a|=|b|$. Thus the similarity class of every quaternion $a$ contains a pair of complex conjugates with absolute-value $|a|$ and real part equal to $\Re( a)$.  Let $a$ be similar to $re^{i \theta}$, $\theta \in (-\pi, \pi]$. We shall adopt the convention of calling $|\theta|$ as the \emph{argument} of $a$ and will denote it by $\arg(a)$.

\subsection{Quaternionic Matrices}
Let ${\rm M}(2, \H)$ denote the group of all $2 \times 2$ quaternionic matrices.  
For $M=\begin{pmatrix}a&b\\c&d\end{pmatrix}\in{\rm M{(2, \H)}}$, define the `quaternionic determinant' of $M$ by 
 $$\det M=|ad-aca^{-1}b|.$$

\medskip \begin{theorem} \cite{kel}, \cite{p} Let $M=\begin{pmatrix}a&b\\c&d\end{pmatrix}\in{\rm M{(2, \H)}}$ be such that $\det M\neq 0$. Then $M$ is invertible and 
$$M^{-1}=\begin{pmatrix} d\sptilde&- b\sptilde\\-c\sptilde&a\sptilde\end{pmatrix}, where $$
\begin{align*}
d\sptilde&=l_{11}^{-1}d, & c\sptilde&=l_{21}^{-1}c,& b\sptilde&=l_{12}^{-1}b,& a\sptilde&=l_{22}^{-1}a; 
\end{align*}
\begin{align*}
l_{11}&=da-dbd^{-1}c & l_{12}&=bdb^{-1}a-bc \\
l_{21}&=cac^{-1}d-cb &  l_{22}&=ad-aca^{-1}b.\end{align*}
\end{theorem}

\medskip  Let
$$\SL(2, \H)=\bigg\{\begin{pmatrix}a&b\\c&d\end{pmatrix}\in {\rm  M}(2, \H):\det{\begin{pmatrix}a&b\\c&d\end{pmatrix}}       =|ad-aca^{-1}b|=1\bigg \}.$$
The group $\SL(2, \H)$ acts by the orientation-preserving isometries of the hyperbolic $5$-space $\h^5$, see \cite{p} for more details.  We identify the extended quaternionic line $\widehat \H=\H \cup \{\infty\}$ to the conformal boundary $\s^4$ of the hyperbolic $5$-space. The group $\S$ acts on $\widehat \H$ by M\"obius transformations:
$$\begin{pmatrix}a&b\\c&d\end{pmatrix}: Z \mapsto (aZ+b)(cZ+d)^{-1}.$$
The action is extended over $\h^5$ by Poincar\'e extensions. Under this action, the group of orientation-preserving isometries of $\h^5$ is ${\rm PSL}(2, \H)=\SL(2, \H)/\{+I, -I\}$. However, often we will not distinguish between an isometry of $\h^5$ and its linear representation in $\SL(2, \H)$. 

\subsection{Classification of isometries}   Every isometry of $\h^5$ has a fixed point on the closure of the hyperbolic space $\overline{\h}^5$ and this gives us the usual classification of elliptic, parabolic and  hyperbolic (or loxodromic) elements in the isometry group. Further, it follows from the  Lefschetz fixed point theorem that every isometry has a fixed point on the conformal boundary. Up to conjugacy, we can take that fixed point to be $\infty$. It follows that every element in $\S$ is conjugate to an upper-triangular matrix. For more details of the classification and algebraic criteria to detect them,  see \cite{cao2}, \cite{kg}, \cite{ps}, also see \cite{fo}.

\subsection{J{\o}rgensen inequality}  The following result is a J{\o}rgensen type inequality for two-generator subgroups of $\SL(2, \H)$ when one of the generators is either elliptic or hyperbolic.

\begin{theorem} \cite{gm} \label{jgi}
Let $ S=\begin{pmatrix}a&b\\c&d\end{pmatrix}$ and
$T=\begin{pmatrix}{\lambda}&0\\0&{\mu}\end{pmatrix}$, $\lambda$ is not similar to $\mu$,  generate a discrete non-elementary subgroup  $\langle S, T \rangle$ of $\S$. Then
 $$\{(\Re\lambda -\Re\mu)^2 +(|\Im \lambda|+|\Im\mu|)^2\}(1+|bc|)\geq 1.$$
\end{theorem}
This gives the following.
\begin{cor}{\rm (\cite{gm}, \cite{kel})}\label{jel}
Let $ S=\begin{pmatrix}a&b\\c&d\end{pmatrix}$ and $T=\begin{pmatrix}{\lambda}&0\\0&{\mu}\end{pmatrix}\in\S$, $\lambda$ is not similar to $\mu$,  generate a discrete non-elementary subgroup $\langle S, T \rangle$  of $\S$. Then
 $$2(\cosh{\tau}-\cos(\alpha+{\beta}))(1+|bc|)\geq 1,$$ where $\alpha=arg(\lambda),~{\beta}=arg(\mu)$, $\tau=2 \log |\lambda|$.
\end{cor}
Observe that with the above expression of $\tau$, we have that $2 \cosh \tau=|\lambda|^2 + |\lambda|^{-2}$.  When one of the generators is  a translation, we have the following result.
\begin{cor}\label{wat} {\rm (\cite{wat}, \cite{kel})}
If  $S=\begin{pmatrix} a&b\\ c&d \end{pmatrix},~T=\begin{pmatrix} 1& \lambda\\ 0& 1\end{pmatrix}$ generate a non-elementary discrete subgroup in $\S$, then $|c| .|\lambda| \geq 1$.
\end{cor}
\subsection{Limit Sets}
Let $L(G)$ be the limit set of a subgroup $G$ of $M(n)$,
see \cite{rat} for basic properties of limit sets.  The limit set $L(G)$ is a closed $G$-invariant subset of $\s^n$. The group $G$ is elementary if $L(G)$ is finite. If $G$ is elementary, $L(G)$ consists of at most two points. If $G$ is non-elementary, then $L(G)$ is an infinite set and every non-empty, closed $G$-invariant subset of $\s^n$ contains $L(G)$.
We note the following lemma, for a proof see \cite[Chapter 12]{rat}.
\begin{lemma}\label{nel} Let $G$ be a subgroup of $M(n)$. 
Let $a \in \partial \h^{n+1}$ be a fixed point of  a non-elliptic element of  $G$.  Then $a$ is a limit point of $G$.
\end{lemma}
Let $F$ be the set of fixed points of all non-elliptic elements of $G$. The above lemma implies that $F$ is $G$-invariant. Further if $G$ is non-elementary, then $F$ contains at least three points. We will use these facts while proving the theorems. Another crucial result to be used in the next section is the following. 
\begin{theorem} \cite[Corollary 4.5.1]{cg}  \label{cgt}
Let $G$ be a subgroup of $\S$ that does not leave invariant a point in $\overline \h^5$ or a proper totally geodesic submanifold in $\h^5$ which is invariant under $G$. Then $G$ is either discrete or dense in $\S$.
\end{theorem}

\section{Discreteness using a Test Map}\label{ms}
\subsection{Proof of \thmref{thm1}}
By hypothesis, $G$ is a Zariski-dense subgroup of $\S$. Therefore, $G$ is non-
elementary. In the sequel we suppose that $G$ is not discrete and derive contradictions
when considering the cases (1)--(3) in the statement of the theorem. 

Suppose  $G$ is not discrete.  Then $G$ is a dense subgroup of $\S$.  It is a well-known fact,  eg.  see \cite{yang4},  that the set of all hyperbolic elements is open in $\SL(2, \H)$. Hence we may choose a hyperbolic  element $g =\begin{pmatrix} a& b\\ c& d\end{pmatrix}$ in $G$ such that it fixes a point other than $0, \infty$.

Let $z_0 \neq 0, \infty$ be a fixed point of $g$.  Consider the element $h = \begin{pmatrix} z_0^{-1}& -1\\ 0& z_0\end{pmatrix}$. It is easy to see that $h^{-1}=\begin{pmatrix} z_0& 1\\ 0& z_0^{-1}\end{pmatrix}$.  Note that $h(z_0)=0$. Since $G$ is dense in $\S$, so there exists a  sequence $\{h_n\} \subseteq G$ such that $h_n \to h$.  We can choose $h_n$   such that $h_n(z_0) \neq 0 \neq h_m(z_0)$ for large $n$, $m$.

\medskip (1) \medskip  Suppose $f$ is $2$-rotatory elliptic.  We can assume, up to conjugacy that,
$$f = \begin{pmatrix} \lambda& 0\\ 0& \mu\end{pmatrix}, ~ \lambda, \mu \in \C, $$
  $|\lambda|=|\mu|=1$,  $\lambda$ is not similar to $\mu$. Further assume  $0 < \at(f)= \arg \lambda + \arg \mu <\frac{\pi}{3}$. Let
 $\arg \lambda= \alpha, \; \arg \mu= {\beta}$.

Let $h_n gh_n^{-1}=\begin{pmatrix} a_n & b_n \\c_n & d_n \end{pmatrix}$.
 By hypothesis, each two generator subgroup $\langle f, h_n gh_n^{-1}\rangle$ is  discrete.  For large $n$,  it follows from \lemref{nel} that $\langle f, h_n g h_n^{-1} \rangle$ has at least three limit points $0$, $\infty$ and $h_n(z_0)$,  and hence, it is non-elementary.   By \thmref{jgi},  for sufficiently large $n$,
 $$2(1- \cos (\alpha + {\beta}))(1+ |b_nc_n|) \geq 1.$$
Now note that
\begin{eqnarray*}
hgh^{-1} &=&\begin{pmatrix} z_0^{-1}& -1\\ 0& z_0\end{pmatrix} \begin{pmatrix} a& b\\ c& d\end{pmatrix} \begin{pmatrix} z_0& 1\\ 0& z_0^{-1}\end{pmatrix}\\
         &=& \begin{pmatrix} z_0^{-1}az_0& z_0^{-1}a+z_0^{-1}bz_0^{-1}-c-dz_0^{-1}\\ z_0cz_0& z_0c+z_0dz_0^{-1}\end{pmatrix}.
\end{eqnarray*}
Since $z_0$ is a fixed point of $g$, we have
\begin{eqnarray*}& &  (az_0+b)(cz_0+d)^{-1}= z_0\\
&\hbox{that is, }&  (z_0^{-1}a+z_0^{-1}bz_0^{-1}-c-dz_0^{-1})z_0cz_0 =0.
\end{eqnarray*}
Since $0<\alpha+{\beta}<\frac{\pi}{3}$, this implies,
\begin{eqnarray*}
2(1-\cos (\alpha+{\beta}))(1+|( z_0^{-1}a+z_0^{-1}bz_0^{-1}-c-dz_0^{-1})z_0cz_0|)&=& 2(1-\cos (\alpha+{\beta}))<1.
\end{eqnarray*}
By \thmref{jgi}, this contradiction completes the proof of (1).

\medskip (2) Let $f$ be hyperbolic. Using the hypothesis, we can assume up to conjugacy that
$$f = \begin{pmatrix} \lambda& 0\\ 0& \mu\end{pmatrix}, \; |\lambda| \neq |\mu|,\; |\lambda \mu|=1, \; \arg \lambda=\alpha, \; \arg \mu ={\beta}, \; 2 \cos(\alpha +{\beta}) > |\lambda|^2+|\mu|^2-1.$$

Let $h_n gh_n^{-1}=\begin{pmatrix} a_n & b_n \\c_n & d_n \end{pmatrix}$.
 By hypothesis and using   \corref{jel}, we have for sufficiently large $n$,
\begin{equation}\label{e1} 2(\cosh \tau- \cos (\alpha + {\beta}))(1+ |b_nc_n|) \geq 1,\end{equation}
where $\tau=2 \log |\lambda|$.
But we have, \begin{eqnarray*}
hgh^{-1} &=&\begin{pmatrix} z_0^{-1}& -1\\ 0& z_0\end{pmatrix} \begin{pmatrix} a& b\\ c& d\end{pmatrix} \begin{pmatrix} z_0& 1\\ 0& z_0^{-1}\end{pmatrix}\\
         &=& \begin{pmatrix} z_0^{-1}az_0& z_0^{-1}a+z_0^{-1}bz_0^{-1}-c-dz_0^{-1}\\ z_0cz_0& z_0c+z_0dz_0^{-1}\end{pmatrix}.
\end{eqnarray*}
Note that $(z_0^{-1}a+z_0^{-1}bz_0^{-1}-c-dz_0^{-1})z_0cz_0 =0$. It follows that
\begin{eqnarray*}
2(\cosh \tau -\cos (\alpha+{\beta}))(1+|( z_0^{-1}a+z_0^{-1}bz_0^{-1}-c-dz_0^{-1})z_0cz_0|)& & \\
= 2(\cosh \tau -\cos (\alpha+{\beta})).
\end{eqnarray*}
 Since, $2 \cos(\alpha +{\beta}) > |\lambda|^2+|\mu|^2-1$, this implies
$$2(\cosh \tau -\cos (\alpha+{\beta}))<1.$$
This is a contradiction to \eqnref{e1}. Hence part (2) of the theorem follows.

\medskip (3)
 Consider the parabolic element $u= \begin{pmatrix} 1& 0\\ -{z_0}^{-1}  & 1\end{pmatrix}$. Note that $u(0)=0$. It is easy to see that $u^{-1}=\begin{pmatrix} 1& 0\\ {z_0}^{-1}& 1\end{pmatrix}$.  Since $G$ is dense in $\S$,  there exists a distinct sequence $\{g_n\} \subseteq G$ such that $g_n \to u$. We may choose $g_n$  such that for large $ n$, $g_n(z_0) \neq \infty$,  and hence,  having  $\langle f, g_n g g_n^{-1}\rangle$  non-elementary.  By hypothesis, these groups are all discrete. Hence,    by \corref{wat},  $$|c_n|.|\mu| \geq 1,$$
where $g_n g g_n^{-1}=\begin{pmatrix} a_n & b_n \\ c_n & d_n \end{pmatrix}$.
By computations we see that 
\begin{eqnarray*}
ugu^{-1} &=&\begin{pmatrix} 1& 0\\ -{z_0}^{-1}& 1\end{pmatrix} \begin{pmatrix} a& b\\ c& d\end{pmatrix} \begin{pmatrix} 1& 0\\ {z_0}^{-1}& 1\end{pmatrix}\\
         &=& \begin{pmatrix} a+b{z_0}^{-1}& b\\ -{z_0}^{-1}(a+b{z_0}^{-1})+(c+d{z_0}^{-1})& -{z_0}^{-1}b+d\end{pmatrix}.
\end{eqnarray*}
Since $z_0$ is a fixed point of $g$, so we have 
$$c_{\infty}= -{z_0}^{-1}(a+b{z_0}^{-1})+(c+d{z_0}^{-1}) =0.$$
 Since $|\mu| \leq 1$, this implies,
\begin{eqnarray*}
|c_n| \geq \frac{1}{|\mu|} &\geq &1.
\end{eqnarray*}
But we see that $c_n \to c_{\infty}=0$ as $n \to \infty, $ which gives a contradiction. This proves (3).  This completes the proof.

\subsection{Proof of \thmref{thm2}}

By similar arguments as used at the beginning of the proof of \thmref{thm1}, we can choose $h_n$   such that $h_n(z_0) \neq 0 \neq h_m(z_0)$ for large $n$, $m$. Let $h_n gh_n^{-1}=\begin{pmatrix} a_n & b_n \\c_n & d_n \end{pmatrix}$. 
\begin{enumerate}
  \item For all $n$, consider 
  \begin{eqnarray*}
    L_n &=& h_ngh_n^{-1}fh_ng^{-1}h_n^{-1} \\
     &=&  \begin{pmatrix} a_n & b_n \\ c_n & d_n \end{pmatrix} \begin{pmatrix} \lambda & 0\\ 0 & \mu \end{pmatrix} \begin{pmatrix}d\sptilde_n&-b\sptilde_n\\-c\sptilde_n&a\sptilde_n\end{pmatrix}\\
     &=& \begin{pmatrix}a_n\lambda d\sptilde_n-b_n\mu c\sptilde_n&-a_n\lambda b\sptilde_n+b_n\mu a\sptilde_n\\
c_n\lambda d\sptilde_n-d_n\mu c\sptilde_n&-c_n\lambda b\sptilde_n+d_n\mu a\sptilde_n\end{pmatrix}\\
     &=& \begin{pmatrix}
           A_n & B_n \\
           C_n & D_n
         \end{pmatrix}. 
  \end{eqnarray*}
As $n \to \infty$, let $L_n \to L_{\infty}$, where $$L_{\infty}=hgh^{-1}fhg^{-1}h^{-1}= \begin{pmatrix}
                                       A_{\infty} & B_{\infty} \\
                                       C_{\infty} & D_{\infty}
                                     \end{pmatrix}.$$
Now we see that,
  \begin{eqnarray*}
    |B_nC_n| & \leq & |a_nb_nc_nd_n||\lambda-a_n^{-1}b_n\mu a\sptilde_n{b\sptilde_n}^{-1}||\lambda-c_n^{-1}d_n\mu c\sptilde_n{d\sptilde_n}^{-1}| \\
     &=& \{(\Re\lambda -\Re\mu)^2 +(|\Im \lambda|+|\Im\mu|)^2\}(1+|b_nc_n|)|b_nc_n|.   \end{eqnarray*}
Let $$\begin{pmatrix}
               a_0 & b_0 \\
               c_0 & d_0
             \end{pmatrix}=hgh^{-1}=\begin{pmatrix} z_0^{-1}az_0& z_0^{-1}a+z_0^{-1}bz_0^{-1}-c-dz_0^{-1}\\ z_0cz_0& z_0c+z_0dz_0^{-1}\end{pmatrix}.$$ 
Since $z_0$ is a fixed point of $g$, we have seen that 
$$ (z_0^{-1}a+z_0^{-1}bz_0^{-1}-c-dz_0^{-1})z_0cz_0 =0,$$
which shows that $b_0c_0=0.$\\
By a similar calculations above in the case $L_n,$ we see that
\begin{eqnarray*}|B_{\infty}C_{\infty}| &\leq& \{(\Re\lambda -\Re\mu)^2 +(|\Im \lambda|+|\Im\mu|)^2\}(1+|b_0c_0|)|b_0c_0|=0, \end{eqnarray*}
and therefore we have $B_{\infty}C_{\infty}=0$. 
This shows that $B_nC_n \to 0.$
Now we see that by hypothesis, each two generator subgroup $\langle f, L_n\rangle$ is  discrete and  non-elementary.  So by \thmref{jgi}, 
 \begin{equation}\label{e2} 2(1- \cos (\alpha + {\beta}))(1+ |B_nC_n|) \geq 1.\end{equation}
Since $0<\alpha+{\beta}<\frac{\pi}{3}$, this implies for sufficiently large $n$,
\begin{eqnarray*}
2(1-\cos (\alpha+{\beta}))(1+|B_nC_n|)&=& 2(1-\cos (\alpha+{\beta}))<1.
\end{eqnarray*}
This is a contradiction to \eqnref{e2} which completes the proof of (1).\\

\medskip  \item For this part the proof follows from similar calculations as in the proof of (1) and the fact that \begin{eqnarray*}
2(\cosh \tau -\cos (\alpha+{\beta}))(1+|B_{\infty}C_{\infty}|)& & \\
= 2(\cosh \tau -\cos (\alpha+{\beta})).
\end{eqnarray*}
 Since, $2 \cos(\alpha +{\beta}) > |\lambda|^2+|\mu|^2-1$, this implies
$$2(\cosh \tau -\cos (\alpha+{\beta}))<1.$$
This leads to a contradiction. Hence part (2) of the theorem follows.\\
 \medskip  \item Consider the parabolic element $h= \begin{pmatrix} 1& 0\\ -{z_0}^{-1}  & 1\end{pmatrix}$. Note that $h(0)=0$. It is easy to see that $h^{-1}=\begin{pmatrix} 1& 0\\ {z_0}^{-1}& 1\end{pmatrix}$.  Since $G$ is dense in $\S$,  there exists a  sequence $\{h_n\} \subseteq G$ such that $h_n \to h$. We may choose distinct $h_n$  such that for large $ n$, $h_n(z_0) \neq \infty$.

Let
\begin{eqnarray*} L_n &=& h_n g h_n^{-1}fh_ng^{-1}h_n^{-1}\\
&=&  \begin{pmatrix} a_n & b_n \\ c_n & d_n \end{pmatrix} \begin{pmatrix} 1 & \mu\\ 0 & 1 \end{pmatrix} \begin{pmatrix}d\sptilde_n&-b\sptilde_n\\-c\sptilde_n&a\sptilde_n\end{pmatrix}\\
     &=& \begin{pmatrix}a_n d\sptilde_n-a_n \mu c\sptilde_n-b_n c\sptilde_n & -a_n\mu a\sptilde_n\\
-c_n\mu c\sptilde_n &-c_n b\sptilde_n+c_n\mu a\sptilde_n+d_n a\sptilde_n\end{pmatrix}\\
     &=& \begin{pmatrix}
           A_n & B_n \\
           C_n & D_n
         \end{pmatrix},\; \text{say}.
\end{eqnarray*}
Now as $n \to \infty$,  $L_n \to L_{\infty}$, where 
\begin{eqnarray*}
  L_{\infty} &=& hgh^{-1}fhg^{-1}h^{-1} \\
   &=& \begin{pmatrix}
         A_{\infty} & B_{\infty} \\
         C_{\infty} & D_{\infty}
       \end{pmatrix},\; \text{say}.
\end{eqnarray*}
It is clear that for large values of $n$, $\langle f, L_n\rangle$  are non-elementary and by hypothesis, these groups are also discrete. Hence,    by \corref{wat},  $|C_n|.|\mu| \geq 1$. Let
$$
  hgh^{-1} =         \begin{pmatrix} a+b{z_0}^{-1}& b\\ -{z_0}^{-1}(a+b{z_0}^{-1})+(c+d{z_0}^{-1})& -{z_0}^{-1}b+d\end{pmatrix}=
\begin{pmatrix}
               a_0 & b_0 \\
               c_0 & d_0
             \end{pmatrix}.$$ We have seen that since $z_0$ is a fixed point of $g$, so
$$c_0= -{z_0}^{-1}(a+b{z_0}^{-1})+(c+d{z_0}^{-1}) =0.
$$
Thus it follows that $C_{\infty}=0$. So,  $C_n \to 0$, as $n \to \infty$.  Since $|\mu| \leq 1$, this implies,
\begin{eqnarray*}
|C_n| \geq \frac{1}{|\mu|} &\geq &1,
\end{eqnarray*}
\end{enumerate} 
 which leads to a contradiction. This completes the proof.

\medskip \begin{ack}
We thank the referee for comments and suggestions. We are also grateful to John Parker for comments on a first draft of this paper. 

\medskip Part of this work was carried out when \hbox{Gongopadhyay}  was visiting the UNSW Sydney supported by the Indo-Australia EMCR Fellowship of the Indian National Science Academy (INSA). Gongopadhyay thanks UNSW for hospitality and the INSA for the fellowship during the visit. Gongopadhyay also thanks Jadavpur University, Kolkata for hospitality where this work was initiated. \end{ack}

\end{document}